
\def\cite{} 
\def\berg{13}
\def\flye{12}
\def\dawson{6}
\def\dickins{2}
\def\euler{9}
\def\frolow{7}
\def\hubers{10}
\def\jaenisch{11}
\def\jello{4}
\def\jelli{5}
\def\jw{8}
\def\lucas{1}
\def\nashw{3}
\def\jr{14}
\def\juengi{15}
\def\juengii{16}

\magnification\magstephalf

\raggedbottom
\baselineskip14pt
\parskip3pt
\font\sc=cmcsc10 

\def\proof{\noindent{\bf Proof.}\quad}
\def\bib{\par\noindent\hangindent 25pt}
\def\biba{\par\parindent 25pt\hangindent 25pt}
\def\pfbox
  {\hbox{\hskip 3pt\lower0pt\vbox{\hrule
  \hbox to 7pt{\vrule height 7pt\hfill\vrule}
  \hrule}}\hskip3pt}
\def\\#1#2{\rm\scriptscriptstyle #1#2\atop
           \rm\scriptscriptstyle #2#1}
\def\longdash{\raise.5ex\hbox to 1.3em{\hrulefill}}

\def\Abar{\overline A}
\def\Bbar{\overline B}
\def\Cbar{\overline C}
\def\Dbar{\overline D}
\def\Ebar{\overline E}
\def\Fbar{\overline F}
\def\Gbar{\overline G}
\def\Hbar{\overline H}
\def\Kbar{\overline K}
\def\Xbar{\overline X}
\def\Obar{\overline O}
\def\ldt{\mathrel{\,.\,.\,}}

\centerline{\bf Leaper Graphs}
\smallskip
\centerline{Donald E. Knuth}
\centerline{Computer Science Department}
\centerline{Stanford University}
\bigskip

\noindent
An $\{r,s\}$-leaper 
[\cite\lucas, p.~130;
\cite\dickins, p.~30; 
\cite\nashw]
is a generalized knight that can jump from $(x,y)$ to $(x\pm r,y\pm s)$ or
$(x\pm s,y\pm r)$ on a rectangular grid. The graph of an $\{r,s\}$-leaper on an
$m\times n$ board is the set of $mn$~vertices $(x,y)$ for $0\leq x<m$ and
$0\leq y<n$, with an edge between vertices that are one $\{r,s\}$-leaper move
apart. We call $x$ the {\it rank\/} and $y$ the {\it file\/} of board position
$(x,y)$.
George~P. Jelliss 
[\cite\jello,
\cite\jelli]
raised several interesting questions about these graphs, and established some
of their fundamental properties. The purpose of this paper is to characterize
when the graphs are connected, for arbitrary~$r$ and~$s$, and to determine the
smallest boards with Hamiltonian circuits when $s=r+1$ or $r=1$.

\proclaim
Theorem 1. The graph of an $\{r,s\}$-leaper on an $m\times n$ board, when
$2\leq m\leq n$ and $1\leq r\leq s$, is connected if and only if the following
three conditions hold: (i)~$r+s$ is relatively prime to $r-s$; (ii)~$n\geq 2s$;
(iii)~$m\geq r+s$.

\proof
Condition (i) is necessary because any common divisor~$d$ of $r+s$ and $r-s$
will be a divisor of $x+y$ for any vertex $(x,y)$ reachable from $(0,0)$; any
leaper move changes the sum of coordinates by $\pm(r+s)$ or $\pm(s-r)$. If
$d>1$, vertex $(0,1)$ would therefore be disconnected from $(0,0)$.

Condition (ii) is necessary because the ``middle'' point $(\lfloor m/2\rfloor,
\lfloor n/2\rfloor)$ would otherwise be isolated: If $m\leq n<2s$ we have
$\lfloor m/2\rfloor-s<0$, $\lfloor m/2\rfloor+s\geq m$, $\lfloor
n/2\rfloor-s<0$, and $\lfloor n/2\rfloor+s\geq n$, so there is no place to
leap from there.

To show that condition (iii) is necessary we show first that the
$\{r,s\}$-leaper graph on any board with $m=r+s-1$ and $n=\infty$ has no path
from $(0,0)$ to $(0,1)$. For this purpose we construct a special path through
points $(x_k,y_k)$ as follows:
$$\eqalign{(x_0,y_0)&=(r-1,0)\,;\cr
\noalign{\smallskip}
(x_{k+1},y_{k+1})&=\cases{(x_k+r,y_k+s)\,,&if $x_k<s-1\,;$\cr
\noalign{\smallskip}
(x_k-s,y_k+r)\,,&if $x_k\geq s\,.$\cr}
\cr}$$
The path terminates when $x_k=s-1$. It is not difficult to see that termination
will occur when $k=r+s-2$, because the sequence $x_0,x_1,\ldots$ runs first
through all values $\{0,1,\ldots,r+s-2\}$ that are congruent to $-1$
modulo~$r$, then all values congruent to $-1-s$, then all values congruent to
$-1-2s$, etc. Since $s$ is relatively prime to~$r$, all residues will occur
before we finally reach $-1-(r-1)s$ modulo~$r$, which is the class of values
congruent to $s-1$. Now if $(x,y)$ is any point reachable from $(r-1,0)$ on
this infinite graph, there is a unique value of~$k$ such that $x=x_k$. And the
difference $y-y_k$ must be an even number, because all $\{r,s\}$-leaper moves
preserve this condition. Therefore the point $(x,y+1)$ is not reachable from
$(r-1,0)$.

Condition (iii) is therefore necessary. Vertex $(0,1)$ will surely be
disconnected from $(0,0)$ on any $m\times n$ board
with $m<r+s$ if it is disconnected
from $(0,0)$ on a $(r+s-1)\times\infty$ board.

Finally we show that conditions (i), (ii), and (iii) are in fact sufficient for
connectivity. We need only prove that the graph is connected when $m=r+s$ and
$n=2s$, because the connectivity of an $m\times n$ board obviously implies
connectivity for $(m+1)\times n$ and $m\times(n+1)$ boards when $m>1$. The
proof is somewhat delicate, because the graph is, in some sense, ``just
barely'' connected.

Let $m=r+s$ and $n=2s$, and let $t$ be any number in the range $0\leq t<n$.
Define a path on the $m\times n$ board by the rules
$$\eqalign{(x_0,y_0)&=(0,t)\,;\cr
\noalign{\smallskip}
(x_{k+1},y_{k+1})&=\cases{(x_k+r,y_k\pm s)\,,&if $x_k<s\,;$\cr
\noalign{\smallskip}
(x_k-s,y_k\pm r)\,,&if $x_k\geq s\,.$\cr}
\cr}$$
The sign of $\pm s$ is uniquely determined by the condition $0\leq y_{k+1}<2s$;
the sign of $\pm r$ may or may not be forced by this condition, and we can use
any desired convention when a choice is possible. The path will reach a point
$(x_k,y_k)$ of the form $(0,u)$ when $k=r+s$, after doing $s$~moves by $(r,\pm
s)$ and $r$~moves by $(-s,\pm r)$. Therefore $(0,u)$ is reachable from $(0,t)$;
we want to use this information to establish connectivity of the graph.

Consider first the case $t=0$; we will choose the signs so that $y_k$ is
either~0 or~$r$ or~$s$ or~$s+r$ for all~$k$. This produces a sequence of
files~$y_k^{(0)}$. Similarly, when $t=r$,~$s$, or~$s+r$, we can keep $y_k$ in
the set $\{0,r,s,s+r\}$, and this defines sequences $y_k^{(r)}$, $y_k^{(s)}$,
$y_k^{(r+s)}$. Now notice that after $k$~steps we have $x_k=ar-bs$ for some
uniquely determined~$a$ and~$b$, where $k=a+b$, and the value of~$y_k$ depends
only on the parity of~$a$ and~$b$:

$$\vcenter{\halign{\hfil#\hfil\qquad%
&\hfil#\hfill\qquad%
&$\hfil#\hfil$\qquad%
&$\hfil#\hfil$\qquad%
&$\hfil#\hfil$\qquad%
&$\hfil#\hfil$\cr
$a$&$b$&y_k^{(0)}&y_k^{(r)}&y_k^{(s)}&y_k^{(r+s)}\cr
\noalign{\smallskip}
even&even&0&r&s&r+s\cr
odd&even&s&r+s&0&r\cr
odd&odd&r+s&s&r&0\cr
even&odd&r&0&r+s&s\cr}}$$

After $r+s$ steps we have $x_{r+s}=0$ and $a=s$, $b=r$. 
Therefore if $r$ is even and $s$ is odd, we have found paths from $(0,0)$ to
$(0,s)$ and from $(0,r)$ to $(0,r+s)$. Traversing these paths twice shows that
$(x,0)$ and $(x,s)$ are reachable from $(0,0)$ whenever $x$ is even; similarly,
$(x,r)$ and $(x,r+s)$ are reachable from $(0,r)$ whenever $x$ is even. In the
other case, when $r$ is odd and $s$ is even, the construction proves that
$(x,0)$ and $(x,r)$ are reachable from $(0,0)$ whenever $x$ is even, while
$(x,s)$ and $(x,r+s)$ are similarly reachable from $(0,s)$. The same argument
establishes a more general principle, which can be formulated as follows:

\proclaim
Lemma 1. Let $t$ be any value such that $0\leq t<t+r<s$. Then $(x,t)$ is
reachable from $(z,t)$ and $(x,t+r)$ is reachable from $(z,t+r)$ whenever
$x-z$ is even. 

\noindent
The proof consists of forming sequences $y_k^{(t)}$, $y_k^{(t+r)}$,
$y_k^{(t+s)}$, and $y_k^{(t+r+s)}$ as before.

Our next step in proving Theorem~1 is to establish a mild form of connectivity:

\proclaim
Lemma 2. Every cell on the $(r+s,2s)$ board is reachable from some cell in
file~0. That is, for all $(x,y)$ there is a $z$ such that $(x,y)$ is connected
to $(z,0)$.

\noindent
Let us say that file $y$ is {\it accessible\/} if its cells are all reachable
from file~0. Let $d=s-r$. To prove that all files are accessible, we will start
at file~0 and, whenever file~$y$ is accessible, we will increase~$y$
by~$d$ if
$y<r$, or decrease~$y$ by~$r$ if $y\geq r$. This procedure will prove
accessibility for all $y<s$. It will then be
obvious that all files $y\geq s$ are accessible.

If $r\leq y<s$ and file $y$ is accessible, we can easily show that file $y-r$
is accessible by using Lemma~1. For if $(x,y-r)$ is any cell in file $y-r$,
there is a path from $(x,y-r)$ to $(0,y-r)$ or to 
$(1,y-r)$, and we can go in one
step from $(z,y-r)$ to $(z+s,y)$ whenever $z<r$.
If $r=z=1$, we need three steps: $(1,y-1)$, $(0,y+s-1)$, $(s,y+s)$, $(s-1,y)$.

The other case is slightly more complicated. Suppose file $y$ is accessible and
$0\leq y<r$. Let $(x,y+d)$ be any cell on file $y+d$. If $x<r$, we can go to
$(x+s,y+d+r)$ to $(x+d,y)$ in two steps, and a similar two-step path applies if
$x\geq s$. So the only problematic case arises when $r\leq x<s$. In such a case
we can follow a zigzag path
$$(x,y+d),\,(x-r,y+d+s),\,(x+d,y+2d),\,\ldots,\,\bigl(x+kd,y+(k+1)d\bigr)$$
until first reaching $y+(k+1)d\geq r$. When this occurs, we have $kd\leq
y+kd<r$, so $x+kd<s+kd<r+s$ is a legitimate rank. Now we can use Lemma~1 to
connect $\bigl(x+kd,y+(k+1)d\bigr)$ to either $\bigl(0,y+(k+1)d\bigr)$ or
$\bigl(1,y+(k+1)d\bigr)$. And we can take another zigzag path from
$\bigl(0,y+(k+1)d\bigr)$ back to file~$y$ as desired:
$$\bigl(0,y+(k+1)d\bigr),\,\bigl(s,y+(k+1)d+r\bigr),\,(d,y+kd),\,\ldots,\,
(kd+s,y+d+r),\,\bigl((k+1)d,y\bigr)\,.$$
Adding $(1,0)$ to each point of this path will also connect
$\bigl(1,y+(k+1)d\bigr)$ to file~$y$, unless $kd=r-1$.

Therefore our proof of Lemma 2 hinges on being able to find a path in the
exceptional case $kd=r-1$. This case can arise only when $y=0$. Hence we must
find a path from $(x,s-1)$ to file~0, for some odd integer~$x$, whenever the
parameters~$r$ and~$s$ satisfy the special conditions $r=1+kd$, $s=1+(k+1)d$.
(Such leaper graphs
 exist whenever $k\geq 0$ and $d$ is odd.) An examination of
small cases reveals a strategy that works in general: The path begins
$$\medmuskip=2mu\eqalign{&(d,s-1),\,(s,2s-1)\,,\cr
&(0,s-1+d),\,(s,s-1+d-r),\,(d,s-1+2d),\,\ldots,\,\bigl(kd,s-1+(k+1)d\bigr)
 =(z-s,2s-s)\,,\cr
&(z,s-2+d),\,(z-s,s-2+d-r),\,(z-d,s-2+2d),\,\ldots,\,
\bigl(z-kd,s-2+(k+1)d\bigr) =(s,2s-3)\,,\cr}$$
where $z=r+s-1$ is the number of the last rank;
the idea is to repeat the $(2k+1)$-step staircase subpaths
$d-1$ times, until reaching $(s,2s-d)$.
Since $2s-d=s+r$, this point is two easy steps from
file~0; the proof of Lemma~2 is complete.

Lemmas 1 and 2 together show that the graph has at most two connected
components, because each vertex is connected either to cell $(0,0)$ or cell
$(1,0)$ of file~0. Furthermore, the construction in the proof of Lemma~2 shows
that $(x,y)$ is connected to $(0,0)$ if $x+f(y)$ is even, to $(1,0)$ if
$x+f(y)$ is odd, where $f(y)$ is a certain parity function associated with
file~$y$. This follows because the proof connects $(x,y)$ to a cell congruent
modulo~2 to $(x+d,y+d)$ when we increase~$y$ by~$d$, to $(x+s,y-r)$ when we
decrease~$y$ by~$r$, and to $(x+r,y+s)$ when we increase~$y$ by~$s$. Thus we
may take $f(y+d)=f(y)+d$, $f(y-r)=f(y)+s$, and $f(y+s)=f(y)+r$.

But the full cycle of changes in $y$ involves $r$ increases by~$d$ and $d$
decreases by~$r$; so we get back to $y=0$ with a parity value of $rd+ds$, which
is odd. Therefore $(0,0)$ is connected to $(1,0)$! The proof of Theorem~1 is
complete. \ \pfbox

\medskip
The next natural question to ask about leaper graphs is whether or not they are
Hamiltonian. Indeed, this is an especially appropriate question, because the
whole idea of Hamiltonian circuits first arose in connection with knight's
tours---which are Hamiltonian paths of the $\{1,2\}$-leaper graph on an
$8\times 8$ board. Knight's tours have fascinated people for more than 1000
years, yet their secrets have not yet been fully unlocked.

Dawson [\cite\dawson]
showed that $\{1,2k\}$-leapers have Hamiltonian paths from corner to corner of
a $(2k+1)\times 4k$ board. Therefore the smallest connected graph already has a
Hamiltonian path when $r=1$. However, it turns out that Hamiltonian circuits
require larger boards, of size $3\times 10$ or $5\times 6$ when $k=1$, or
$9\times 10$ when $k=2$. 
Jelliss [\cite\jelli]
derived necessary conditions for the existence of Hamiltonian circuits using
$\{r,r+1\}$-leapers, and conjectured that the board of smallest area in this
special case has size $(2r+1)\times (6r+4)$. He proved this conjecture when
$r\leq 3$.

The following theorem gives further support to Jelliss's conjecture, because it
shows that a board with smallest dimension greater than $2r+1$ must have an
area at least $(4r+2)^2$; this exceeds $(2r+1)(6r+4)$.

\proclaim
Theorem 2. If\/ $r>2$ and the graph of an $\{r,r+1\}$-leaper on an\/
 $m\times n$ board has a Hamiltonian circuit, and if\/
 $2r+1<m\leq n$, then\/ $m\geq4r+2$.

\proof
Let $s=r+1$. All vertices $(x,y)$ with $0\leq x,y<r$ are adjacent only to
$(x+r,y+s)$ and $(x+s,y+r)$. Therefore a Hamiltonian circuit must include both
of these edges. Similarly, the edges from $(x,y)$ to $(x-r,y+s)$ and
$(x-s,y+r)$ are forced when $m-r\leq x<m$ and $0\leq y<r$.

Therefore, if $2r+1<m<3r+1$, a ``snag''
[\cite\jelli]
occurs at vertex $(r,s)$: Any Hamiltonian circuit must lead from this vertex 
directly to
$(0,0)$, but also to $(2r,0)$ and $(r+s,1)$, because $m-r\leq 2r<r+s=2r+1<m$.
Only two of these three compulsory moves are possible.

Similarly, if $3r+1<m<4r+1$, there is a snag at $(2r,s)$. This vertex must
connect to $(r-1,1)$, $(3r,0)$, and $(3r+1,1)$. 

The case $m=3r+1$ is impossible if $r>1$, because vertex $(2r-1)$ must connect
to $(r-2,1)$, $(3r-1,0)$, and $(3r,1)$.

Suppose finally that $m=4r+1$. Vertex $(r-1,2)$ has just two neighbors,
$(2r-1,r+3)$ and $(2r,r+2)$, because $r>2$. Similarly, $(3r+1,2)$ is adjacent
only to $(2r+1,r+3)$ and $(2r,r+2)$. Therefore $(2r,r+2)$ must not connect to
$(3r,1)$. But now $(3r,1)$ has only two remaining options, namely $(4r,r+2)$
and $(2r-1,r+1)$. This makes a snag at $(2r-1,r+1)$, which must also link to
$(r-1,0)$ and $(r-2,1)$. \ \pfbox

\medskip
On the other hand, Jelliss's conjecture does turn out to be false for all $r>3$,
because $\{r,r+1\}$-leapers continue to acquire new problems on narrow boards
as $r$ grows:

\proclaim
Theorem 3. If $r>3$ and the graph of an $\{r,r+1\}$-leaper on a $(2r+1)\times
n$ board has a Hamiltonian circuit, then $n\geq r^2+5r+2$ if $r$ is odd,
$r^2+6r+4$ if $r$ is even.

\proof
Let $s=r+1$, and suppose first that we have a board of size
$(2r+1)\times\infty$. Certain edges are forced to be in any Hamiltonian
circuit, because some vertices have degree~2. We will see that such edges,
in turn, can force other connections. 

Each vertex $(x,y)$ has at most four neighbors. If the $x$~coordinate
represents vertical position (rank) and the $y$~coordinate represents
horizontal position 
(file) as in matrix notation, two of these neighbors lie to the
``left,'' namely $(x\pm r,y-s)$ and $(x\pm s,y-r)$, and two lie to the
``right,'' $(x\pm r,y+s)$ and $(x\pm s,y+r)$. (If $x\neq r$, there is one
choice of sign for $x\pm s$ and $x\pm r$; if $x=r$, neither choice works for
$x\pm s$, but both choices are valid for $x\pm r$.)

It will be convenient to use a two-dimensional representation of the files: The
notation $[a,b]$ will refer to $as-b$, for $a\geq 1$ and $1\leq b\leq s$. In
terms of this notation, the neighbors of all vertices in file $[a,b]$ belong to
files $[a-1,b]$, $[a-1,b-1]$, $[a+1,b]$, and $[a+1,b+1]$. (Appropriate
adjustments to these formulas are made when $a=1$ or $b=1$ or $b=s$.)

The files are also classified into various types:

\medskip
\bib
Type $R$: All vertices must link to both right neighbors. (For example, file
$[1,s]=0$ is type~$R$ because its vertices have no left neighbors. In fact,
file $[1,b]$ is type~$R$ for $2\leq b\leq s$.)

\bib
Type $R'$: The middle vertex, in rank $r$, must link to both right neighbors.
(File $[1,1]$ is type~$R'$, because vertex $(r,r)$ has no left neighbors; the
other vertices in this file have one  neighbor to the left.)

\bib
Type $L$: All vertices must link to both left neighbors. (File $[2,s]$ is
type~$L$, because both left neighbors of its vertices are in files of type~$R$.
In fact, file $[2,b]$ is type~$L$ for $3\leq b\leq s$.)

\bib
Type $L'$: The middle vertex must link to both left neighbors. (File $[2,2]$ is
type~$L'$.)

\bib
Type $l$: All vertices except the middle must link to at least one left
neighbor. (File $[1,1]$ is type~$l$, because for example vertex $(0,r)$ has a
left neighbor in file~0, which is type~$R$.)

\bib
Type $l'$: The extreme vertices, in ranks 0 and~$2r$, must link to at least one
left neighbor. $\bigl($File $[2,1]=2s-1$ is type~$l'$, because $(0,2s-1)$ and
$(2r,2s-1)$ are the only neighbors of $(r,r)$.$\bigr)$

\medskip
Additional types $r$ and $r'$ arise when $n$ is finite, because the right
boundary has properties like that of the left.

When file $[a,b]$ is type~$R$, its right neighbors $[a+1,b]$ and $[a+1,b+1]$
must be either type~$L$ or~$l$. When $[a,b]$ is type~$R'$, $[a+1,b]$ is
either~$L$ or~$l'$. And when files $[a,b]$ and $[a,b+1]$ both have type~$L$,
file $[a+1,b]$ must have type~$R$; its left neighbors cannot be used. Such
arguments inductively establish the following facts for all $a\leq s$:

\medskip
\biba
$[a,b]$ is type~$R$ when $a$ is odd, type~$L$ when $a$ is even, for $a<b\leq
s$; 

\biba
$[a,a]$ is type $R'$ when $a$ is odd, type $L'$ when $a$ is even;

\biba
$[a,1]$ is type $l$ when $a$ is odd and $a<s$;

\biba
$[a,a]$ is type $l$ when $a$ is even;

\biba
$[a,a-1]$ is type $l'$ when $a$ is even.

\medskip
\noindent
Furthermore $[s+1,s]$ is type $l'$ when $r$ is even. A file can simultaneously
have two non-conflicting types; for example, $[1,1]$ is both $R'$ and~$l$,
while $[2,2]$ is both~$L'$ and~$l$. At the other extreme, files with no forced
links have no special type.

To prove the theorem, we will show first that all files having special types
when $n=\infty$ must be present when $n$ is finite. Theorem~1 tells us that
$n\geq 2s$; we want to prove that $n>s^2$. Suppose $n=ks-d$, where $0\leq d<s$
and $k\leq s$. This right boundary introduces complementary constraints; let 
$$\overline{[a,b]}=n-1-[a,b]=[k-a,d+1-b]$$
be the file that corresponds to $[a,b]$ when left and right are interchanged.
Then file $\overline{[a,b]}$ has the type we have derived for $[a,b]$ but with
the interchange of~$L$ and~$R$, $l$~and~$r$. The value of $n=ks-d$ must be
chosen so that complementary types do not conflict.

First, $k$ must be even. For if $k$ is odd, file $(k-1)s=[k,s]$ is type~$R$
or~$R'$; but it is also file $\overline{[1,d+1]}$, which is type~$L$ or~$L'$.

Second, $d$ must be even. For the leaper graphs are always bipartite---each
leap links a ``black'' square with a ``red'' square, as on a chessboard---and a
bipartite graph cannot be Hamiltonian when it has an odd number of vertices.

Third, $d$ must be zero. Otherwise file $[k,s]$ is type~$L$ or~$l$, while it is
also $\overline{[1,1]}$, which is type~$r$. Since $L$ is incompatible with~$r$,
we must have $k=s$; files $[1,1]$ and $s,s]$ will then be of type~$lr$. This,
however, forces a ``short circuit,''
$$\eqalign{&(0,r),\,(r,r+s),\,(2r,r)\,,\cr
&\qquad
(r-1,0),\,(2r-1,r+1),\,(r-2,1),\,\ldots,\,(0,r-1),\,(r,2r),\,(2r,r-1)\,,\cr 
&\qquad
(r-1,2r-1),\,(2r-1,r-2),\,(r-2,2r-2),\,\ldots,
\,(1,r+1),\,(r+1,0),\,(0,r)\,.\cr}$$

Fourth, we must have $k>d$. Otherwise file $[k-1,k-1]$, which is type~$R'$,
would link to the nonexistent file $[k,k-1]\geq n$.

Fifth, a contradiction arises even when $k$ and $d$ are even and $0<d<k\leq s$.
Suppose, for example, that $k=10$ and $d=6$. All neighbors of vertices in
Hamiltonian circuits are then forced except for certain vertices in 24~files:
$$\vcenter{\halign{%
$\hfil#\hfil$\quad%
&$\hfil#\hfil$\quad%
&$\hfil#\hfil$\quad&$\hfil#\hfil$\quad&$\hfil#\hfil$\quad&$\hfil#\hfil$\quad%
&$\hfil#\hfil$\quad&$\hfil#\hfil$\quad&$\hfil#\hfil$\quad&$\hfil#\hfil$\quad%
&$\hfil#\hfil$\cr
[1,1]&\longdash&[2,2]&\longdash&[3,3]&\longdash&[4,4]&\longdash%
&[5,5]&\longdash&[6,6]\cr
\vert&&\vert&&\vert&&\vert&&\vert&&\vert\cr
[2,1]&\longdash&[3,2]&\longdash&[4,3]&\longdash&[5,4]&\longdash&[6,5]%
&\longdash&[7,6]\cr
\vert&&\vert&&\vert&&\vert&&\vert&&\vert\cr
[3,1]&\longdash&[4,2]&\longdash&[5,3]&\longdash%
&[6,4]&\longdash&[7,5]&\longdash&[8,6]\cr
\vert&&\vert&&\vert&&\vert&&\vert&&\vert\cr
[4,1]&\longdash&[5,2]&\longdash&[6,3]&\longdash%
&[7,4]&\longdash&[8,5]&\longdash&[9,6]\cr}}$$
This follows because all files $[a,b]$ with $a<b$ are type~$R$ or~$L$; and all
files $\overline{[a,b]}$ with $a<b$ are type~$L$ or~$R$. We have
$\overline{[a,b]}=[10-a,7-b]$, so the files in the array can also be
represented in the dual form
$$\vcenter{\halign{%
$\hfil#\hfil$\quad%
&$\hfil#\hfil$\quad%
&$\hfil#\hfil$\quad&$\hfil#\hfil$\quad&$\hfil#\hfil$\quad&$\hfil#\hfil$\quad%
&$\hfil#\hfil$\quad&$\hfil#\hfil$\quad&$\hfil#\hfil$\quad&$\hfil#\hfil$\quad%
&$\hfil#\hfil$\cr
\overline{[9,6]}&\longdash&\overline{[8,5]}&\longdash&\overline{[7,4]}%
&\longdash&\overline{[6,3]}&\longdash%
&\overline{[5,2]}&\longdash&\overline{[4,1]}\cr
\vert&&\vert&&\vert&&\vert&&\vert&&\vert\cr
\overline{[8,6]}&\longdash&\overline{[7,5]}&\longdash&\overline{[6,4]}%
&\longdash&\overline{[5,3]}&\longdash%
&\overline{[4,2]}&\longdash&\overline{[3,1]}\cr
\vert&&\vert&&\vert&&\vert&&\vert&&\vert\cr
\overline{[7,6]}&\longdash&\overline{[6,5]}&\longdash&\overline{[5,4]}%
&\longdash&\overline{[4,3]}&\longdash&\overline{[3,2]}%
&\longdash&\overline{[2,1]}\cr
\vert&&\vert&&\vert&&\vert&&\vert&&\vert\cr
\overline{[6,6]}&\longdash&\overline{[5,5]}&\longdash&\overline{[4,4]}%
&\longdash&\overline{[3,3]}&\longdash&\overline{[2,2]}%
&\longdash&\overline{[1,1]}\cr}}$$
All links between files are indicated by horizontal and vertical lines in these
arrays.

Let us say that a vertex is even or odd according as it belongs to file $[a,b]$
where $a$ is even or odd. All links go between even vertices and odd vertices.
The edges of a Hamiltonian circuit that have not been forced by our arguments
so far must therefore touch the same number of even vertices as odd vertices.
We will obtain a contradiction by showing that the odd vertices have more
unspecified neighbors than the even vertices~do.

In our example, the 24 files are classified as follows:
$$\vcenter{\halign{%
$\hfil#\hfil$\quad%
&$\hfil#\hfil$\quad%
&$\hfil#\hfil$\quad&$\hfil#\hfil$\quad&$\hfil#\hfil$\quad&$\hfil#\hfil$\quad%
&$\hfil#\hfil$\quad&$\hfil#\hfil$\quad&$\hfil#\hfil$\quad&$\hfil#\hfil$\quad%
&$\hfil#\hfil$\cr
R'l&\longdash&L'l&\longdash&R'&\longdash&L'l&\longdash&R'&\longdash&L'l\cr
\vert&&\vert&&\vert&&\vert&&\vert&&\vert\cr
l'&\longdash&\emptyset&\longdash&l'&\longdash&\emptyset&\longdash&l'%
&\longdash&r\cr
\vert&&\vert&&\vert&&\vert&&\vert&&\vert\cr
l&\longdash&r'&\longdash&\emptyset&\longdash&r'&\longdash%
&\emptyset&\longdash&r'\cr
\vert&&\vert&&\vert&&\vert&&\vert&&\vert\cr
R'r&\longdash&L'&\longdash&R'r&\longdash&L'&\longdash&R'r&\longdash&L'r\cr}}$$
An unconstrained file (indicated here by $\emptyset$) has $2r+1$ vertices with
a total of $4r+2$ unspecified neighbors. Making it type~$l$ or~$r$
specifies~$2r$ of these; making it type~$l'$, $r'$, $L'$, or~$R'$ specifies~2.
The total number of specified neighbors in odd files is
$4(2r)+6(2)$; the total in even files  is $6(2r)+12(2)$. Therefore
the odd vertices have an excess of unspecified neighbors; there aren't enough
``slots'' available to specify them all.

For general $k$ and $d$, there are $(k-d)d$ files with partially unspecified
vertices. There will be $d/2$ odd files of type~$L'$, and $d/2$ of type~$R'$:
the same holds for even files. There are $(k-d)/2$ odd files of type~$l$, and
$(k-d)/2$ of type~$r$; there are $d/2$ even files of type~$l$, and $d/2$ of
type~$r$. There are no odd files of types~$l'$ or~$r'$; there are $d/2$ even
files of each of those types. Therefore the total number of unspecified
vertices in even files will balance the total in odd files if and only if
$$(k-d)(2r)=d(2r+2)\,.$$
The smallest solution to this equation, when $d>0$ and $k-d>0$, is $k-d=r+1$,
$d=r$. But then $k=2r+1$ exceeds~$s$.

We know therefore that $n>s^2$. Hence all the files of special types, from
$[1,s]=0$ to $[s+1,s]=s^2$, are present. Further considerations depend on
whether $r$ is even or odd.

Suppose $r$ is even. Then the arguments above can be used also in the case
$k=s+1$, except that when $d=0$ we find that file $[1,1]$ is of type~$lr'$
instead of~$lr$. Contradictions are obtained exactly as before; hence
$n>(s+1)s$. 

We cannot have $n=(s+2)s-d$ for $d=3,5,\ldots,r-1$, because file
$\overline{[s,s]}$, which is type~$L'$, would be the same as file $[3,d+1]$,
which is type~$R$. The case $d=1$ is also impossible, because file
$\overline{[s+1,s]}$, which has type~$r'$, is file $[2,2]$, which has
type~$L'l$; then vertices $(0,2r)$ and $(2r,2r)$ would both link to $(r,r-1)$
and $(3r,r+1)$, forming a short circuit. Therefore $n>(s+2)s$.

We can show $n\geq (s+3)s$ by letting $n=(s+3)s-d$ and using the parity
argument above with $k=s+3$. When $d=2$ the number of odd files of type~$l$
(and of type~$r$) will be one less than before, because $[s,1]$ is not
type~$l$, but this just makes the lack of balance even worse. Furthermore, when
$d=0$ and $r>4$, file $\overline{[r-1,1]}=[6,s]$ has conflicting types~$r$
and~$L$.

When $r=4$ and $n=40=(s+3)s$, file $\overline{[s+1,s]}=[3,1]$ has types~$r'$
and~$l$; this forces file $[2,1]$ to be type~$L'$, and a loop is forced from
$(r,2r+1)$ to $(2r,r)$ to \dots to $(0,r)$ to $(r,2r+1)$. 

And when $n=(s+4)s-d$ for $d=5,7,\ldots,r-1$, file $\overline{[s,s]}=[5,d+1]$
has conflicting types~$L'$ and~$R$. 
If $d=3$, file $\overline{[s+1,s]}=[4,4]$ is~$r'$ and~$l$, forcing a loop from
$(r,3r-1)$ to $(0,4r)$ to $(r,5r+1)$ to $(r,4r)$ to $(r,3r-1)$. 

Finally, suppose $r$ is odd. We cannot have $n=(s+1)s-d$ for $d=2,4,\ldots,r-1$
because file $\overline{[s,s]}=[2,d+1]$ has conflicting types~$R'r$ and~$L$.
The case $d=0$ is also impossible, because file $\overline{[s,s-1]}=[2,2]$
would be of types~$r'$ and~$L'l$, forcing a short circuit as before. Next, if
$n=(s+2)s-d$, we use the parity argument with $k=s+2$ when $d>0$. If $d=0$ and
$r>3$, file $\overline{[r,1]}=[4,s]$ has conflicting types~$r$ and~$L$. And
when $n=(s+3)s-d$ for $d=4,6,\ldots,r-1$, file $\overline{[s,s]}=[4,d+1]$ has
conflicting types~$R'r$ and~$L$. \ \pfbox

\medskip
Theorems 2 and 3 have proved that certain leaper graphs fail to be Hamiltonian.
Let us now strike a happier note by constructing infinitely many Hamiltonian
circuits, on the boards of smallest area not ruled out by those theorems.

\proclaim
Theorem 4. The graph of an $\{r,r+1\}$-leaper on a $(4r+2)\times (4r+2)$ board
is Hamiltonian.

\proof
Let us call the leaper moves {\sc ne}, {\sc nw}, {\sc en},
{\sc es}, {\sc se}, {\sc sw}, {\sc ws}, and~{\sc wn},
where $(\hbox{\sc n},\hbox{\sc e},\hbox{\sc s},\hbox{\sc w})$
 stand respectively for North, East, South, West, and the first letter
indicates the direction of longest leap. North is the direction of
decreasing~$x$; East is the direction of increasing~$y$. Thus {\sc ne}
 is a leap from $(x,y)$ to $(x-s,y+r)$; {\sc es} goes to $(x+r,y+s)$. 

We construct first a highly symmetric set of leaper moves in which every vertex
has degree~2, illustrated here in the case $r=4$:

$$\vcenter{\halign{$\hfil#\hfil$\ \ %
&$\hfil#\hfil$\ \ %
&$\hfil#\hfil$\ \ 
&$\hfil#\hfil$\ \ %
&$\hfil#\hfil$\ \ %
&$\hfil#\hfil$\ \ %
&$\hfil#\hfil$\ \ %
&$\hfil#\hfil$\ \ %
&$\hfil#\hfil$\ \ %
&$\hfil#\hfil$\ \ %
&$\hfil#\hfil$\ \ %
&$\hfil#\hfil$\ \ %
&$\hfil#\hfil$\ \ %
&$\hfil#\hfil$\ \ %
&$\hfil#\hfil$\ \ %
&$\hfil#\hfil$\ \ %
&$\hfil#\hfil$\ \ %
&$\hfil#\hfil$\cr
\\SE&\\SE&\\SE&\\SE&D_1&\\SW&\\SW&\\SW&A_4&\Abar_4&\\SE&\\SE&\\SE&\Dbar_1%
&\\SW&\\SW&\\SW&\\SW\cr
\noalign{\vskip4pt}
\\SE&\\SE&\\SE&\\SE&D_2&\\SW&\\SW&\\SW&A_3&\Abar_3&\\SE&\\SE&\\SE&\Dbar_2%
&\\SW&\\SW&\\SW&\\SW\cr
\noalign{\vskip4pt}
\\SE&\\SE&\\SE&\\SE&D_3&\\SW&\\SW&\\SW&A_2&\Abar_2&\\SE&\\SE&\\SE&\Dbar_3%
&\\SW&\\SW&\\SW&\\SW\cr
\noalign{\vskip4pt}
\\SE&\\SE&\\SE&\\SE&D_4&\\SW&\\SW&\\SW&A_1&\Abar_1&\\SE&\\SE&\\SE&\Dbar_4%
&\\SW&\\SW&\\SW&\\SW\cr
\noalign{\vskip4pt}
B_1&B_2&B_3&B_4&F_4&E_4&E_3&E_2&E_1&\Ebar_1&\Ebar_2&\Ebar_3&\Ebar_4&\Fbar_4%
&\Bbar_4&\Bbar_3&\Bbar_2&\Bbar_1\cr
\noalign{\vskip4pt}
\\NE&\\NE&\\NE&\\NE&F_3&\\NW&\\NW&\\NW&\\NW&\\NE&\\NE&\\NE&\\NE&\Fbar_3%
&\\NW&\\NW&\\NW&\\NW\cr
\noalign{\vskip4pt}
\\NE&\\NE&\\NE&\\NE&F_2&\\NW&\\NW&\\NW&\\NW&\\NE&\\NE&\\NE&\\NE&\Fbar_2%
&\\NW&\\NW&\\NW&\\NW\cr
\noalign{\vskip4pt}
\\NE&\\NE&\\NE&\\NE&F_1&\\NW&\\NW&\\NW&\\NW&\\NE&\\NE&\\NE&\\NE&\Fbar_1%
&\\NW&\\NW&\\NW&\\NW\cr
\noalign{\vskip4pt}
X&\\NE&\\NE&\\NE&C_0&C_1&C_2&C_3&C_4&\Cbar_4&\Cbar_3&\Cbar_2&\Cbar_1&\Cbar_0%
&\\NW&\\NW&\\NW&\Xbar\cr
\noalign{\vskip4pt}
X'&\\SE&\\SE&\\SE&C'_0&C'_1&C'_2&C'_3&C'_4&\Cbar'_4&\Cbar'_3%
&\Cbar'_2&\Cbar'_1&\Cbar'_0%
&\\SW&\\SW&\\SW&\Xbar'\cr
\noalign{\vskip4pt}
\\SE&\\SE&\\SE&\\SE&F'_1&\\SW&\\SW&\\SW&\\SW&\\SE&\\SE&\\SE&\\SE&\Fbar'_1%
&\\SW&\\SW&\\SW&\\SW\cr
\noalign{\vskip4pt}
\\SE&\\SE&\\SE&\\SE&F'_2&\\SW&\\SW&\\SW&\\SW&\\SE&\\SE&\\SE&\\SE&\Fbar'_2%
&\\SW&\\SW&\\SW&\\SW\cr
\noalign{\vskip4pt}
\\SE&\\SE&\\SE&\\SE&F'_3&\\SW&\\SW&\\SW&\\SW&\\SE&\\SE&\\SE&\\SE&\Fbar'_3%
&\\SW&\\SW&\\SW&\\SW\cr
\noalign{\vskip4pt}
B'_1&B'_2&B'_3&B'_4&F'_4&E'_4&E'_3&E'_2&E'_1&\Ebar'_1&\Ebar'_2%
&\Ebar'_3&\Ebar'_4&\Fbar'_4%
&\Bbar'_4&\Bbar'_3&\Bbar'_2&\Bbar'_1\cr
\noalign{\vskip4pt}
\\NE&\\NE&\\NE&\\NE&D'_4&\\NW&\\NW&\\NW&A'_1&\Abar'_1&\\NE&\\NE&\\NE&\Dbar'_4%
&\\NW&\\NW&\\NW&\\NW\cr
\noalign{\vskip4pt}
\\NE&\\NE&\\NE&\\NE&D'_3&\\NW&\\NW&\\NW&A'_2&\Abar'_2&\\NE&\\NE&\\NE&\Dbar'_3%
&\\NW&\\NW&\\NW&\\NW\cr
\noalign{\vskip4pt}
\\NE&\\NE&\\NE&\\NE&D'_2&\\NW&\\NW&\\NW&A'_3&\Abar'_3&\\NE&\\NE&\\NE&\Dbar'_2%
&\\NW&\\NW&\\NW&\\NW\cr
\noalign{\vskip4pt}
\\NE&\\NE&\\NE&\\NE&D'_1&\\NW&\\NW&\\NW&A'_4&\Abar'_4&\\NE&\\NE&\\NE&\Dbar'_1%
&\\NW&\\NW&\\NW&\\NW\cr
}}$$

\noindent
Some vertices have been given names $A_j$, $B_j$, etc.; these provide important
connecting links. The other, nameless vertices occur in blocks where all links
have the same pair of directions; for example, every vertex called 
$\\SE$ is joined to its {\sc se} and {\sc es} neighbors. These blocks 
effectively serve as
parallel mirrors that provide staircase paths rising or falling at
$45^{\circ}$~angles.

The named vertices are limited as follows:

\smallskip
\halign{\qquad\qquad#\hfil\cr
$A_j$ goes {\sc ws}, {\sc es};\cr
$B_j$ goes {\sc se}, {\sc en};\cr
$C_j$ goes {\sc nw}, {\sc sw}, except $C_r$ goes {\sc nw}, {\sc se};\cr
$D_j$ goes {\sc se}, {\sc sw};\cr
$E_j$ goes {\sc wn}, {\sc ws}, except $E_r$ goes {\sc wn}, {\sc se};\cr
$F_j$ goes {\sc nw}, {\sc en}, except $F_r$ goes {\sc sw}, {\sc en};\cr
$X$ goes {\sc ne}, {\sc se}.\cr}

\noindent
The directions for complemented vertices like $\Abar_j$ are the same but with
{\sc e} and~{\sc w} interchanged; the directions for primed vertices
like~$A'_j$ are the same but with {\sc n} and~{\sc s} interchanged. Thus, for
example, $B'_j$ goes {\sc ne}, {\sc es}; and $\Fbar'_r$ goes {\sc ne},
{\sc ws}\null. 
It is easy to verify that these pairs of directions are consistent. 

We can now deduce the connections between named vertices, following paths
through unnamed ones. The nearest named neighbors of~$A_j$ are $B_j$
and~$\Fbar_j$. The nearest to~$B_j$ are~$A_j$ and~$C'_{j-1}$. The nearest
to~$C_j$ are~$B'_{j+1}$ (or $\Ebar'_r$ when $j=r$) and~$D_j$ (or~$E_1$
when $j=0$). The nearest to~$D_j$ are~$C_j$ and~$E_j$ (or~$X$ when $j=r$). The
nearest to~$E_j$ are~$D_j$ (or~$\Cbar'_r$ when $j=r$) and $F_{j-1}$
(or~$C_0$ when $j=1$). The nearest to~$F_j$ are~$\Abar_j$ and~$E_{j+1}$
(or~$X'$ when $j=r$). And the nearest to~$X$ are~$D_r$ and~$F'_r$.

Each chain of edges therefore falls into a pattern that depends on the value
of~$r\bmod 4$. If $r=4k$ we have
$$\eqalign{X\ldt D_r&\ldt C_r\ldt \Ebar'_r\ldt\Fbar'_{r-1}
\ldt A'_{r-1}\ldt B'_{r-1}\ldt C_{r-2}\ldt D_{r-2}\cr
&\qquad \ldt E_{r-2}\ldt F_{r-3}\ldt \Abar_{r-3}\ldt\Bbar_{r-3}\ldt
\Cbar'_{r-4}\ldt\Dbar'_{r-4}\ldt\Ebar'_{r-4}\;\cdots\cr
&\qquad \ldt E_2\ldt F_1\ldt\Abar_1\ldt\Bbar_1\ldt\Cbar'_0\ldt\Ebar'_1
\ldt \Dbar'_1\ldt\Cbar'_1\ldt\Bbar_2\ldt\Abar_2\;\cdots\cr
&\qquad\ldt\Abar_{r-2}\ldt F_{r-2}\ldt E_{r-1}\ldt D_{r-1}\ldt C_{r-1}\ldt
B'_r\ldt A'_r\ldt \Fbar'_r\ldt\overline{X}\,.\cr}$$
If $r=4k+2$ the pattern is almost the same except that the middle transition is
complemented and primed:
$$\ldt\Ebar'_2\ldt\Fbar'_1\ldt A'_1\ldt B'_1\ldt C_0\ldt E_1\ldt D_1\ldt
C_1\ldt B'_2\ldt A'_2\;.$$
If $r=4k+1$ the pattern in the middle is
$$\ldt E_3\ldt F_2\ldt\Abar_2\ldt\Bbar_2\ldt\Cbar'_1\ldt\Dbar'_1\ldt\Ebar'_1
\ldt\Cbar'_0\ldt\Bbar_1\ldt\Abar_1\ldt\;\cdots$$
ending with $\Abar_r\ldt F_r\ldt X'$. And if $r-4k+3$, the middle is again
complemented and primed.

Consequently the edges defined above make exactly two circuits altogether.
If $r$ is even, one circuit contains $X$ and~$\overline{X}$, the other
contains~$X'$ and~$\overline{X}'$. If $r$ is odd, the circuits contain
$\{X,X'\}$ and $\{\overline{X},\overline{X}'\}$, respectively.

A small change now joins the circuits together into a single Hamiltonian
circuit. We simply replace the subpaths
$$F'_r,\,X,\,D_r,\,C_r,\,\Ebar'_r\qquad\quad{\rm and}\quad\qquad
\Fbar_r,\,\overline{X}',\,\Dbar'_r,\,\overline{X}',\,\Ebar'_r\,,$$
by
$$F'_r,\,C_r,\,D_r,\,X,\,E_r\qquad\quad{\rm and}\qquad\quad
\Fbar_r,\,\Cbar'_r,\,\Dbar'_r,\,\Xbar',\,\Ebar'_r\,,$$
respectively. \ \pfbox

\medskip
The proof of Theorem 4 leads, for example, to the following $18\times
18$ leaper tour when $r=4$:
$$\vcenter{\halign{%
$\hfil#$\ %
&$\hfil#$\ \ &$\hfil#$\ \ &$\hfil#$\ \ &$\hfil#$\ \ &$\hfil#$\ \ %
&$\hfil#$\ \ &$\hfil#$\ \ &$\hfil#$\ \ &$\hfil#$\ \ &$\hfil#$\ \ %
&$\hfil#$\ \ &$\hfil#$\ \ &$\hfil#$\ \ &$\hfil#$\ \ &$\hfil#$\ \ %
&$\hfil#$\ \ &$\hfil#$\cr
0&272&220&43&53&333&363&104&183&83&4&263&233&153&143&320&372&100\cr
270&222&41&55&212&51&331&365&102&2&265&231&151&312&155&141&322&370\cr
224&38&57&210&277&214&48&328&367&267&228&148&314&377&310&157&138&324\cr
36&60&207&280&386&275&216&46&326&226&146&316&375&286&380&307&160&136\cr
334&362&105&182&84&388&273&218&44&144&318&373&88&184&82&5&262&234\cr
52&332&364&103&1&271&221&42&54&154&142&321&371&101&3&264&232&152\cr
213&50&330&366&268&223&40&56&211&311&156&140&323&368&266&230&150&313\cr
276&215&47&327&225&37&58&208&278&378&308&158&137&325&227&147&315&376\cr
387&274&217&45&35&61&206&281&385&287&381&306&161&135&145&317&374&285\cr
85&174&117&345&335&361&106&181&87&185&81&6&261&235&245&17&74&187\cr
176&115&347&27&125&337&358&108&178&78&8&258&237&25&127&247&15&76\cr
113&350&30&66&168&123&340&356&111&11&256&240&23&68&166&130&250&13\cr
352&32&64&203&301&171&121&342&354&254&242&21&71&201&303&164&132&252\cr
34&62&205&282&384&288&173&118&344&244&18&73&188&284&382&305&162&134\cr
336&360&107&180&86&175&116&346&26&126&246&16&75&186&80&7&260&236\cr
124&338&357&110&177&114&348&28&67&167&128&248&14&77&10&257&238&24\cr
170&122&341&355&112&351&31&65&202&302&165&131&251&12&255&241&22&70\cr
300&172&120&343&353&33&63&204&283&383&304&163&133&253&243&20&72&200\cr
}}$$
(Radix 9 notation is used here so that the near-fourfold symmetry is revealed.)
This tour was found by an exhaustive computer calculation, which determined
that exactly 16~different leaper tours on this board have $180^{\circ}$
symmetry; none have $90^{\circ}$ symmetry. The chosen tour exhibited maximum
symmetry under the circumstances, and fortunately it could be generalized to
arbitrary~$r$. 

Theorem~4 applies in particular when $r=1$. We get the well
known $6\times 6$ knight's tour
$$\vcenter{\halign{\hfil#\quad&\hfil#\quad&\hfil#\quad%
&\hfil#\quad&\hfil#\quad&\hfil#\cr
0&36&13&23&6&30\cr
12&24&38&28&14&22\cr
37&1&35&7&31&5\cr
25&11&27&15&21&17\cr
2&34&8&18&4&32\cr
10&26&3&33&16&20\cr
}}$$
(again in radix 9).

A. H. Frost showed a century ago that the graphs for $\{1,4\}$-leapers and
$\{2,3\}$-leapers are Hamiltonian on a $10\times 10$ board
[\cite\frolow, plate VII].
T. H. Willcocks showed more recently that $\{2,5\}$-leapers and
$\{3,4\}$-leapers are Hamiltonian on a $14\times 14$
[\cite\jw].
Willcocks conjectured that an $\{r,s\}$-leaper has a Hamiltonian circuit on a
$2(r+s)\times 2(r+s)$ board whenever $s-r$ and $s+r$ are relatively prime.
Theorem~4 establishes infinitely many cases of this conjecture, and computer
calculations have verified it whenever $r+s<15$. The computer had to work hard
only in the case $r=5$, $s=8$.

We can verify Willcocks's conjecture also in the other extreme case, when $r=1$:

\proclaim Theorem 5. The graph of a $\{1,2k\}$-leaper on a $(4k+2)\times
(4k+2)$ board is Hamiltonian.

\goodbreak
\proof
Let $s=2k$. We give symbolic names to the vertices as follows, illustrated here
in the case $k=4$:
$$\vcenter{\halign{%
$\hfil#\hfil$\kern.65em
&$\hfil#\hfil$\kern.65em
&$\hfil#\hfil$\kern.65em
&$\hfil#\hfil$\kern.65em
&$\hfil#\hfil$\kern.65em
&$\hfil#\hfil$\kern.65em
&$\hfil#\hfil$\kern.65em
&$\hfil#\hfil$\kern.65em
&$\hfil#\hfil$\kern.65em
&$\hfil#\hfil$\kern.65em
&$\hfil#\hfil$\kern.65em
&$\hfil#\hfil$\kern.65em
&$\hfil#\hfil$\kern.65em
&$\hfil#\hfil$\kern.65em
&$\hfil#\hfil$\kern.65em
&$\hfil#\hfil$\kern.65em
&$\hfil#\hfil$\kern.65em
&$\hfil#\hfil$\kern.65em
&$\hfil#\hfil$\cr
O&.&.&.&.&.&.&H_2&K_4&.&.&.&B_{22}&.&B_{12}&.&A_2&\Obar\cr
\noalign{\smallskip}
K_5&.&.&.&B_{21}&.&B_{11}&.&A_1&.&.&.&.&.&.&H_1&K_3&.\cr
\noalign{\smallskip}
.&G_2&.&G_4&.&G_6&.&G_8&.&.&.&.&.&.&.&.&.&.\cr
\noalign{\smallskip}
.&.&W_{11}^0&X_{11}^0&W_{12}^0&X_{12}^0&W_{13}^0&X_{13}^0&.&.%
&Y_{11}^1&Z_{11}^1&Y_{12}^1&Z_{12}^1&Y_{13}^1&Z_{13}^1&J_{12}&.\cr
\noalign{\smallskip}
.&F_{11}&Y_{11}^0&Z_{11}^0&Y_{12}^0&Z_{12}^0&Y_{13}^0&Z_{13}^0&J_{11}&.%
&W_{11}^1&X_{11}^1&W_{12}^1&X_{12}^1&W_{13}^1&X_{13}^1&.&.\cr
\noalign{\smallskip}
.&.&W_{21}^0&X_{21}^0&W_{22}^0&X_{22}^0&W_{23}^0&X_{23}^0&.&E_{13}%
&Y_{21}^1&Z_{21}^1&Y_{22}^1&Z_{22}^1&Y_{23}^1&Z_{23}^1&I_{13}&.\cr
\noalign{\smallskip}
.&E_{12}&Y_{21}^0&Z_{21}^0&Y_{22}^0&Z_{22}^0&Y_{23}^0&Z_{23}^0&.&.%
&W_{21}^1&X_{21}^1&W_{22}^1&X_{22}^1&W_{23}^1&X_{23}^1&.&.\cr
\noalign{\smallskip}
.&.&.&D_2&.&D_4&.&D_6&.&.&.&.&.&.&.&.&H_4&.\cr
\noalign{\smallskip}
.&.&.&.&.&.&.&.&H_3&.&.&B_{23}&.&B_{13}&.&A_3&.&.\cr
\noalign{\smallskip}
.&K_6&.&B_{26}&.&B_{16}&.&A_4&.&.&.&.&.&.&.&.&.&K_2\cr
\noalign{\smallskip}
G_1&.&G_3&.&G_5&.&G_7&.&.&K_1&.&.&.&.&.&.&.&.\cr
\noalign{\smallskip}
.&.&X_{11}^2&W_{11}^2&X_{12}^2&W_{12}^2&X_{13}^2&W_{13}^2&F_{13}&.%
&Z_{11}^3&Y_{11}^3&Z_{12}^3&Y_{12}^3&Z_{13}^3&Y_{13}^3&.&J_{13}\cr
\noalign{\smallskip}
F_{12}&.&Z_{11}^2&Y_{11}^2&Z_{12}^2&Y_{12}^2&Z_{13}^2&Y_{13}^2&.&.%
&X_{11}^3&W_{11}^3&X_{12}^3&W_{12}^3&X_{13}^3&W_{13}^3&.&.\cr
\noalign{\smallskip}
.&.&X_{21}^2&W_{21}^2&X_{22}^2&W_{22}^2&X_{23}^2&W_{23}^2&.&.%
&Z_{21}^3&Y_{21}^3&Z_{22}^3&Y_{22}^3&Z_{23}^3&Y_{23}^3&.&I_{12}\cr
\noalign{\smallskip}
E_{11}&.&Z_{21}^2&Y_{21}^2&Z_{22}^2&Y_{22}^2&Z_{23}^2&Y_{23}^2&.&I_{11}%
&X_{21}^3&W_{21}^3&X_{22}^3&W_{22}^3&X_{23}^3&W_{23}^3&.&.\cr
\noalign{\smallskip}
.&.&D_1&.&D_3&.&D_5&.&D_7&.&.&.&.&.&.&.&.&H_5\cr
\noalign{\smallskip}
.&.&.&.&.&.&.&.&.&.&C_2&.&B_{24}&.&B_{14}&.&.&.\cr
\noalign{\smallskip}
O'&.&C_1&.&B_{25}&.&B_{15}&.&.&.&.&.&.&.&.&.&.&\Obar'\cr}}$$
Positions marked `.' obtain names by complementation (left-right reflection)
and/or priming (up-down reflection). Thus, for example, the full names of the
positions on the top row are
$$O\ \Abar_2\ C_1'\ \Bbar_{12}\ B'_{25}\ \Bbar_{22}\ B'_{15}\ H_2\ K_4
\ \Kbar_4\ \Hbar_2\ \Bbar'_{15}\ B_{22}\ \Bbar'_{25}\ B_{12}\ \Cbar'_1\ 
A_2\ \Obar\,;$$
the dots help keep the diagram from being more cluttered than it already is.
The names in general are
$$\eqalign{&A_1,\ldots,A_4\,;\ B_{i1},\ldots,B_{i6}\,;\ C_1,\,C_2\,;\ 
D_1,\ldots,D_{s-1}\,;\cr
\noalign{\smallskip}
&E_{p1},\,E_{p2},\,E_{p3}\,;\ F_{q1},\,F_{q2},\,F_{q3}\,;\ 
G_1,\ldots,G_s\,;\ H_1,\ldots,H_5\,;\cr
\noalign{\smallskip}
&I_{p1},\,I_{p2},\,I_{p3}\,;\ J_{q1},\,J_{q2},\,J_{q3}\,;\ 
K_1,\ldots,K_6\,;\ W^a_{ij},\,X^a_{ij},\,Y^a_{ij},\,Z^a_{ij}\,;\cr}$$
here $1\leq a\leq 4$, \ $1\leq i\leq k-2$, \ $1\leq j\leq k-1$, \
$1\leq p<\lceil k/2\rceil$, and $1\leq q<\lfloor k/2\rfloor$.

The names $E_{pl}$ occur two positions to the left of
$Z^a_{(k-2p)1}$, where $a=(2,0,1)$ for
$l=(1,2,3)$ respectively.
The names $F_{ql}$ occur two positions to the left of
$Z^a_{(k-2q-1)1}$, where $a=(0,2,3)$.
The names $I_{pl}$ occur two positions to the right of
$Y^a_{(k-2p)(k-1)}$, where $a=(2,3,1)$ if $k$ is
even, $a=(0,1,3)$ if $k$ is odd;
the names $J_{ql}$ occur two positions to the right of
$Y^a_{(k-2q-1)(k-1)}$, where $a=(0,1,3)$ if $k$ is
even, $a=(2,3,1)$ if $k$ is odd.
The positions of the other names are self-evident. We may assume that $k\geq
3$. 

To show that the graph is Hamiltonian, we will first link the vertices in six
closed circuits, then we will join those circuits together. The basic circuit,
if $k=2l+2$,~is
$$\eqalign{O,\;&A_1,\ldots,A_4,\,
B_{11},\ldots,B_{16},\ldots,B_{(k-2)1},\ldots,B_{(k-2)6},\;
C_1,C_2,\,D_1,\ldots,D_{s-1}\,,\cr
\noalign{\smallskip}
&E_{11},E_{12},E_{13},\;F_{11},F_{12},F_{13},\;
E_{21},\,\ldots,\,F_{l3},\;G_1,\ldots,G_s,\;H_1,\ldots,H_5\,,\cr
\noalign{\smallskip}
&I_{11},I_{12},I_{13},\;J_{11},J_{12},J_{13},\;
I_{21},\,\ldots,\,J_{l3},\;K_1,\ldots,K_6\,,\cr
\noalign{\smallskip}
O',\;&A'_1,\;\ldots,\;K'_6,\;O\,.\cr}$$
On the other hand, if $k=2l+1$, the basic circuit is
$$\eqalign{O,\;&A_1,\ldots,A_4,\,
B_{11},\ldots,B_{16},\ldots,B_{(k-2)1},\ldots,B_{(k-2)6},\;
C_1,C_2,\,D_1,\ldots,D_{s-1}\,,\cr
\noalign{\smallskip}
&E_{11},E_{12},E_{13},\;F_{21},F_{12},F_{13},\;
E_{21},\,\ldots,\,E_{l3},\;H'_5,\ldots,H'_1,\;G'_s,\ldots,G'_1\,,\cr
\noalign{\smallskip}
&I_{11},I_{12},I_{13},\;J_{11},J_{12},J_{13},\;
I_{21},\,\ldots,\,I_{l3},\;K_1,\ldots,K_6\,,\cr
\noalign{\smallskip}
O',\;&A'_1,\;\ldots,\;K'_6,\;O\,.\cr}$$
Another basic circuit is obtained by complementing everything. We also form
$W$, $X$, $Y$, and~$Z$ circuits as follows, when $i=k-2$ and $j=k-1$:
$$\eqalign{&W^0_{11},W^1_{11},W^0_{21},W^1_{21},\ldots,W^0_{i1},W^1_{i1},
\alpha_i,\alpha_{i-1},\ldots,\alpha_1\,,\cr
\noalign{\smallskip}
&\qquad \alpha_l=W^3_{l1},W^1_{l2},W^3_{l2},\ldots,W^1_{lj},W^3_{lj},W^2_{lj},
W^0_{lj},W^2_{li},\ldots,W^0_{l2},W^2_{l1}\,;\cr
\noalign{\medskip}
&X^2_{11},X^3_{11},X^2_{21},X^3_{21},\ldots,X^2_{i1},X^3_{i1},
\beta_i,\beta_{i-1},\ldots,\beta_1\,,\cr
\noalign{\smallskip}
&\qquad \beta_l=X^1_{l1},X^3_{l2},X^1_{l2},\ldots,X^3_{lj},X^1_{lj},X^0_{lj},
X^2_{lj},X^0_{li},\ldots,X^2_{l2},X^0_{l1}\,;\cr
\noalign{\medskip}
&Y^1_{11},Y^0_{11},Y^1_{21},Y^0_{21},\ldots,Y^1_{i1},Y^0_{i1},
\gamma_i,\gamma_{i-1},\ldots,\gamma_1\,,\cr
\noalign{\smallskip}
&\qquad \gamma_l=Y^2_{l1},Y^0_{l2},Y^2_{l2},\ldots,Y^0_{lj},Y^2_{lj},Y^3_{lj},
Y^1_{lj},Y^3_{li},\ldots,Y^1_{l2},Y^3_{l1}\,;\cr
\noalign{\medskip}
&Z^3_{11},Z^2_{11},Z^3_{21},Z^2_{21},\ldots,Z^3_{i1},Z^2_{i1},
\delta_i,\delta_{i-1},\ldots,\delta_1\,,\cr
\noalign{\smallskip}
&\qquad \delta_l=Z^0_{l1},Z^2_{l2},Z^0_{l2},\ldots,Z^2_{lj},Z^0_{lj},Z^1_{lj},
Z^3_{lj},Z^1_{li},\ldots,Z^3_{l2},Z^1_{l1}\,.\cr}$$
Circuits can be spliced together when we have consecutive vertices $(u_1,u_2)$
in one circuit and $(v_1,v_2)$ in another, where $u_1$ is adjacent to~$v_1$ and
$u_2$ is adjacent to~$v_2$. The pairs 
$$\vcenter{\halign{$#$\hfil\qquad&$#$\hfil\cr
(E_{12},E_{13})&(Z^2_{i1},Z^3_{i1})\cr
\noalign{\smallskip}
(E'_{12},E'_{13})&(W^0_{11},W^1_{11})\cr
\noalign{\smallskip}
(\Dbar_2,\Dbar_3)&(Y^2_{ii},Y^0_{ij})\cr
\noalign{\smallskip}
(\Dbar'_2,\Dbar'_3)&(X^0_{1i},X^2_{1j})\cr
\noalign{\smallskip}
(\Gbar'_3,\Gbar'_3)&(Z^2_{ij},Z^0_{ij})\cr}}$$
satisfy this property and suffice to complete the proof. For the first two
pairs hook the $Z$ and~$W$ 
circuits into the basic circuit, and the next two hook the
$Y$ and~$W$ circuits into its complement; the last pair
hooks the $Z$~circuit, which
is now part of the basic circuit, into the complement. \ \pfbox

\medskip
As an example of the construction in Theorem~5, here is the the $22\times 22$
circuit that arises when $s=10$:

\medskip

\halign{&\kern.4em\hfil#\hfil\cr
  0&143&386&139&384&133&378&127&372&457&360&269&172&257& 18&251& 12&245&  6&243&  2&145\cr
361&270&171&258& 17&252& 11&246&  5&242&  1&144&387&140&383&134&377&128&371&458&359&268\cr
 54&467&388&465&390&463&392&461&394&459&396&185&170&187&168&189&166&191&164&241&162&333\cr
397&274&399&196&401&194&403&192&405&238& 55&332& 81&320& 83&318& 85&316& 87&314&355&184\cr
158&331& 80&307& 94&309& 92&311& 90&313&354&275&446&231&412&233&410&235&408&237& 56&471\cr
449&276&445&230&415&228&417&226&419&224&157&472& 79&306& 97&304& 99&302&101&300&353&180\cr
156&327& 78&293&108&295&106&297&104&299&450&179&444&217&426&219&424&221&422&223& 60&473\cr
451&280&443&216&429&214&431&212&433&210& 61&326& 77&292&111&290&113&288&115&286&349&178\cr
152& 35& 76& 37& 74& 39& 72& 41&118&285&348&281&442&203&440&205&438&207&436&209& 62&477\cr
347&174&345& 26&343& 28&341& 30&339& 32&151&478& 65& 48& 67& 46& 69& 44&119&284&455&282\cr
150&483&142&385&138&379&132&373&126&367&456&173&262& 19&256& 13&250&  7&244&  3&146&479\cr
271&362&263& 22&259& 16&253& 10&247&  4&335& 52&141&382&135&376&129&370&123&366&267&358\cr
468& 53&466&389&464&391&462&393&460&395&272&357&186&169&188&167&190&165&240&163&334&161\cr
273&398&197&400&195&402&193&404&239&406&469&160&321& 82&319& 84&317& 86&315& 88&183&356\cr
330&159&322& 95&308& 93&310& 91&312& 89&182&447&198&413&232&411&234&409&236&407&470& 57\cr
277&448&199&414&229&416&227&418&225&420&329& 58&323& 96&305& 98&303&100&301&102&181&352\cr
328&155&324&109&294&107&296&105&298&103&278&351&200&427&218&425&220&423&222&421&474& 59\cr
279&452&201&428&215&430&213&432&211&434&475&154&325&110&291&112&289&114&287&116&177&350\cr
 34&153& 36& 75& 38& 73& 40& 71& 42&117&176&453&202&441&204&439&206&437&208&435&476& 63\cr
175&346& 25&344& 27&342& 29&340& 31&338& 33& 64& 49& 66& 47& 68& 45& 70& 43&120&283&454\cr
482&149& 50&137&380&131&374&125&368&121&364&265& 24&261& 20&255& 14&249&  8&337&480&147\cr
363&264& 23&260& 21&254& 15&248&  9&336&481&148& 51&136&381&130&375&124&369&122&365&266\cr
}

\medskip
Leapers with $r=1$ can in fact tour a slightly smaller board:

\proclaim
Theorem 6. The graph of a $\{1,2k\}$-leaper on a $(4k+1)\times (4k+2)$ board is
Hamiltonian.

\proof
This time the construction is simpler. We may assume that $k\geq 2$, because
Euler
[\euler] constructed a $5\times 6$ knight's tour. The case $k=2$ was solved by
Huber-Stockar [\hubers], whose method can be generalized to all larger
values of~$k$. We assign names as follows, illustrated when $k=4$:

$$\vcenter{\halign{%
$\hfil#\hfil$\kern.65em
&$\hfil#\hfil$\kern.65em
&$\hfil#\hfil$\kern.65em
&$\hfil#\hfil$\kern.65em
&$\hfil#\hfil$\kern.65em
&$\hfil#\hfil$\kern.65em
&$\hfil#\hfil$\kern.65em
&$\hfil#\hfil$\kern.65em
&$\hfil#\hfil$\kern.65em
&$\hfil#\hfil$\kern.65em
&$\hfil#\hfil$\kern.65em
&$\hfil#\hfil$\kern.65em
&$\hfil#\hfil$\kern.65em
&$\hfil#\hfil$\kern.65em
&$\hfil#\hfil$\kern.65em
&$\hfil#\hfil$\kern.65em
&$\hfil#\hfil$\kern.65em
&$\hfil#\hfil$\kern.65em
&$\hfil#\hfil$\cr
.&\Xbar_4&X_1&A_1^3&A_1^4&A_1^5&A_1^6&A_1^7&.&Y_1&.&.&.&.&.&.&.&Z_2\cr
\noalign{\smallskip}
.&Y_6&.&.&.&.&.&.&.&Z_3&X_2&A_2^3&A_2^4&A_2^5&A_2^6&A_2^7&.&Y_2\cr
\noalign{\smallskip}
X_7&.&X_9&A_5^3&A_5^4&A_5^5&A_5^6&A_5^7&.&Y_9&.&.&.&.&.&.&X_5&.\cr
\noalign{\smallskip}
.&Y_{14}&.&.&.&.&.&.&X_6&.&X_{10}&A_6^3&A_6^4&A_6^5&A_6^6&A_6^7&.&Y_{10}\cr
\noalign{\smallskip}
X_{15}&.&X_{17}&A_9^3&A_9^4&A_9^5&A_9^6&A_9^7&.&Y_{17}&.&.&.&.&.&.&X_{13}&.\cr
\noalign{\smallskip}
.&Y_{22}&.&.&.&.&.&.&X_{14}&.&X_{18}&A_{10}^3&A_{10}^4%
&A_{10}^6&A_{10}^6&A_{10}^7&.&Y_{18}\cr
\noalign{\smallskip}
X_{23}&.&X_{25}&A_{13}^3&A_{13}^4&A_{13}^5&A_{13}^6&A_{13}^7%
&.&Y_{25}&.&.&.&.&.&.&X_{21}&.\cr
\noalign{\smallskip}
.&Y_{30}&.&.&.&.&.&.&X_{22}&.&X_{26}&A_{14}^3&A_{14}^4%
&A_{14}^6&A_{14}^6&A_{14}^7&.&Y_{26}\cr
\noalign{\smallskip}
Z_5&.&B_2&B_3&B_4&B_5&B_6&B_7&B_8&\Bbar_8&\Bbar_7&\Bbar_6&\Bbar_5%
&\Bbar_4&\Bbar_3&\Bbar_2&Z_1&.\cr
\noalign{\smallskip}
Y_5&.&Y_7&.&.&.&.&.&.&X_3&A_3^3&A_3^4&A_3^5&A_3^6&A_3^7&.&Y_3&.\cr
\noalign{\smallskip}
.&X_8&A_4^3&A_4^4&A_4^5&A_4^6&A_4^7&.&Y_4&.&Y_8&.&.&.&.&.&.&X_4\cr
\noalign{\smallskip}
Y_{13}&.&Y_{15}&.&.&.&.&.&.&X_{11}&A_7^3&A_7^4&A_7^5&A_7^6&A_7^7&.&Y_{11}&.\cr
\noalign{\smallskip}
.&X_{16}&A_8^3&A_8^4&A_8^5&A_8^6&A_8^7&.&Y_{12}&.&Y_{16}%
&.&.&.&.&.&.&X_{12}\cr
\noalign{\smallskip}
Y_{21}&.&Y_{23}&.&.&.&.&.&.&X_{19}&A_{11}^3&A_{11}^4&A_{11}^5&A_{11}^6%
&A_{11}^7&.&Y_{19}&.\cr
\noalign{\smallskip}
.&X_{24}&A_{12}^3&A_{12}^4&A_{12}^5&A_{12}^6&A_{12}^7%
&.&Y_{20}&.&Y_{24}&.&.&.&.&.&.&X_{20}\cr
\noalign{\smallskip}
Y_{29}&.&Y_{31}&.&.&.&.&.&.&X_{27}&A_{15}^3&A_{15}^4&A_{15}^5&A_{15}^6%
&A_{15}^7&.&Y_{27}&.\cr
\noalign{\smallskip}
.&Z_6&A_{16}^3&A_{16}^4&A_{16}^5&A_{16}^6&A_{16}^7%
&.&Y_{28}&.&Y_{32}&.&.&.&.&.&.&X_{28}\cr}}$$

\noindent
As in the proof of Theorem 5, dots stand for names obtained by complementation
(left-right reflection). The vertex names, in general, are
$$\eqalign{&A_1^j,\ldots,A_{4k}^j\,,\quad\hbox{for}\quad 2<j<2k\,;\quad
B_2,\ldots,B_{2k}\,;\cr
\noalign{\smallskip}
&X_1,\ldots,X_{8k-4}\,;\quad Y_1,\ldots,Y_{8k}\,;\quad Z_1,\ldots,Z_6\,.\cr}$$
Notice that the graph contains paths
$$(B_{j-1}\hbox{ or }B_{j+1}),\,A_1^j,\,\ldots,\,A_{4k}^j,\,
(B_{j-2}\hbox{ or }B_j)$$
for $2<j<2k$, except that vertex $B_{j-2}$ is not present when $j=3$.

Let $\alpha_j$ be the path
$B_{j-2},A_{4k}^j,\ldots,A_1^j,B_{j+1},A_{4k}^{j+1},\ldots,A_1^{j+1}$, when
$2<j<2k$ and $j$ is even. This path $\alpha_j$ can be followed by
vertex~$B_{j+2}$. Therefore we can get from $X_1$ to~$B_{2k-2}$ and to~$B_{2k}$
via the disjoint paths
$$\eqalign{&X_1,\ldots,X_{8k-4},\,Z_1,\ldots,Z_6,\,\alpha_4,\alpha_8,\ldots,
\alpha_{2k-4},\,B_{2k-2}\cr
\noalign{\smallskip}
\hbox{and \quad}&X_1,\,B_3,\,A_{4k}^3,\ldots,A_1^3,\,\alpha_6,\alpha_{10},\ldots,
\alpha_{2k-2},\,B_{2k}\cr}$$
when $k$ is even,
$$\eqalign{&X_1,\ldots,X_{8k-4},\,Z_1,\ldots,Z_6,\,\alpha_4,\alpha_8,\ldots,
\alpha_{2k-2},\,B_{2k}\cr
\noalign{\smallskip}
\hbox{and \quad}&X_1,\,B_3,\,A_{4k}^3,\ldots,A_1^3,\,\alpha_6,\alpha_{10},\ldots,
\alpha_{2k-4},\,B_{2k-2}\cr}$$
when $k$ is odd. Since $Y_1,\ldots,Y_{8k}$ is a path from $B_{2k}$
to~$\Bbar_{2k-2}$, we obtain a path from~$B_{2k-2}$ to~$\Bbar_{2k-2}$ that runs
through all vertices with uncomplemented names. This path plus its complement
is the desired Hamiltonian circuit. \ \pfbox

\medskip
The board in Theorem 6 turns out to be as small as possible.

\proclaim
Theorem 7. \kern-3.8pt
A $\{1,2k\}$-leaper has no Hamiltonian circuit on a board of
area less than $(4k{+}1)(4k{+}2)$.

\proof
Consider the $\{1,2k\}$-leaper graph on an $m\times n$ board with $m\le n$ and
$2k+1\le m\le 4k$. If $m$ is even, we can show that no Hamiltonian circuit
exists by using an argument due to de Jaenisch [\cite\jaenisch, page~46],
Flye Sainte-Marie [\cite\flye],  and
Jelliss [\cite\jw]: Say that vertex
$(x,y)$ is type~$A$ if $y$~is even, $x$~is even, and $x<2k$, or if $y$~is odd,
$x$~is odd, and $x\geq m-2k$; it is of type~$B$ if $y$~is even, $x$~is odd, and
$x<2k$, or if $y$~is odd, $x$~is even, and $x\geq m-2k$; it is type~$C$
otherwise. Type~$A$ vertices are adjacent only to vertices of type~$B$, but
there are no more~$B$'s than~$A$'s. Therefore the only possible circuit
containing all the~$A$'s has the form $ABAB\ldots AB$. But such a circuit
misses all the~$C$'s.

Suppose therefore that $m$ is odd, say $m=2l+1$. Then $l\ge k$ and
$n\ge4k$, by Theorem~1. The links from $(l,y)$ to $(l\pm1,y+2k)$ are
forced when $y<2k$, because $l<2k$. The case $n=4k$ is
impossible by the argument in the previous paragraph, when the ranks
and files of the board are transposed.
Therefore $n>4k$, and a short
circuit from $(l,0)$ to $(l\pm1,2k)$ to $(l,4k)$ is forced unless
$n>6k$. Indeed, if $6k\le n<8k$, a short circuit from $(l,n-6k)$ to
$(l\pm1,n-4k)$ to $(l,n-2k)$ is forced. Consequently we have
$n\ge8k$.  (This argument, in the case $k=2$, was suggested
by Jelliss in a letter to the author.)

If $m\ge2k+3$ we have therefore $mn\ge16k^2+24k>(4k+1)(4k+2)$. And if
$m=2k+1$ and $n\ne8k$ we have $mn\ge(2k+1)(8k+2)=(4k+1)(4k+2)$ because
$n$~must be even.

The remaining case is quite interesting, because we will see that a
Hamiltonian path (but not a circuit) is possible. Let $m=2k+1$ and $n=8k$,
and assume that a Hamiltonian circuit exists. We will write $u\sim v$ if $u$
and~$v$ are adjacent vertices of the circuit. Vertices of degree~2 force the
connections
$$(x,y)\sim(x+1,y+2k)\,,\qquad (x,y)\sim(x-1,y+2k)\,,$$
for $0<x<2k$ and $0\le y<2k$; also at the corners we have
$$(0,0)\sim(1,2k)\,,\quad(0,0)\sim(2k,1)\,,\qquad
  (2k,0)\sim(0,1)\,,\quad(2k,0)\sim(2k-1,2k)\,.$$
Notice that both neighbors of $(x,y)$ in the circuit have now been identified
whenever $1<x<2k-1$ and $y<4k$; by symmetry, the same is true in the right
half of the board, when $y\ge4k$. Our goal is to deduce the behavior of the
circuit on the remaining vertices, which lie in the top two and bottom two
ranks of the board. We will assume that $k>1$, so that these four ranks are
distinct. A similar (and much simpler) argument applies when $k=1$.

\def\ol(#1){\overline{(#1)}}
Let $\ol(x,y)=(2k-x,8k-1-y)$ be the point opposite $(x,y)$ with respect to
the center of the board. Whenever we deduce that $u\sim v$, a symmetric
derivation proves that $\overline u\sim\overline v$; such consequences need
not be stated explicitly.

We must have either $(0,4k-1)\sim(2k,4k-2)$ or
$(0,4k-1)\sim(2k,4k)=\ol(0,4k-1)$, because
$(0,4k-1)\not\sim(1,6k-1)=\ol(2k-1,2k)$; the latter is joined to
$\ol(2k,0)$ and $\ol(2k,2)$. Similarly, either $(2k,4k-1)\sim(0,4k-2)$ or
$(2k,4k-1)\sim\ol(2k,4k-1)$. These choices are not independent. For if
$(0,4k-1)\sim(2k,4k-2)$ and $(2k,4k-1)\sim(0,4k-2)$, there is a short circuit
$$\eqalign{(0,4k-1)\sim(1,2k-1)\sim\cdots&\sim(2k-1,2k-1)\sim(2k,4k-1)\cr
 &\sim(0,4k-2)\sim\cdots\sim(2k,4k-2)\sim(0,4k-1)\,;\cr}$$
likewise the connections $(0,4k-1)\sim\ol(0,4k-1)$ and $(2k,4k-1)\sim
\ol(2k,4k-1)$ force a short circuit
$$(0,4k-1)\sim\cdots\sim(2k,4k-1)\sim\ol(2k,4k-1)\sim\cdots\sim
 \ol(0,4k-1)\sim(0,4k-1)\,.$$
By symmetry we can therefore assume without loss of generality that
$$(0,4k-1)\sim(2k,4k-2)\qquad\hbox{and}\qquad(2k,4k-1)\sim\ol(2k,4k-1)\,.$$
These connections imply also $\ol(0,4k-1)\sim\ol(2k,4k-2)$; we still are
able to claim legitimately below that $\overline u\sim\overline v$ holds
whenever we have deduced that $u\sim v$.

Further detective work establishes $(1,2k+1)\sim(0,1)$, because
$(1,2k+1)\not\sim(0,4k+1)=\ol(2k,4k-2)$. Therefore the circuit contains the
path
$$(2k,1)\sim(0,0)\sim\cdots\sim(2k,0)\sim(0,1)\sim(1,2k+1)\sim(2,1)\sim\cdots
\sim(2k-2,1)\sim(2k-1,2k+1)\,;$$
it follows that $(2k,1)\not\sim(2k-1,2k+1)$. The only possibilities
remaining are $(2k,1)\sim(0,2)$ and $(2k-1,2k+1)\sim(2k,4k+1)=\ol(0,4k-2)$.

Now we can establish, in fact, the relations
$$(1,2k+2j-1)\sim(0,2j-1)\sim(2k,2j-2)\,,\qquad(0,4k-2j+1)\sim(2k,4k-2j)\,,$$
for $1\le j\le k$. They have been verified when $j=1$; suppose we know them
for some $j<k$. Then $(2k,2j)\not\sim(0,2j-1)$, hence $(2k,2j)\sim(0,2j+1)$
and $(2k,2j)\sim(2k-1,2k+2j)$. Hence $(0,4k-2j-1)\not\sim(1,6k-2j-1)=
\ol(2k-1,2k+2j)$; we must have $(0,4k-2j-1)\sim(2k,4k-2j-2)$. This in turn
forces $(1,2k+2j+1)\sim(0,2j+1)$, because $(1,2k+2j+1)$ can't be joined to
$(0,4k+2j+1)=\ol(2k,4k-2j-2)$. The induction on~$j$ is complete, and we have
also proved
$$(2k,2j)\sim(2k-1,2k+2j)\qquad\hbox{for $1\le j<k$}\,.$$

One consequence of our deductions so far is the existence of a rather long
path,
$$\displaylines{\quad
 (0,2)\sim(2k,1)\sim\cdots\sim(2k-1,2k+1)\sim\ol(0,4k-2)\sim\cdots\sim
 \ol(2k,4k-2)\sim\ol(0,4k-1)\hfill\cr
\hfill\sim\cdots\sim\ol(2k,4k-1)\sim(2k,4k-1)\sim\cdots\sim\ol(0,2)\,.
 \quad\cr}$$
We've also found paths from $(0,2k+2j-2)$ to $(2k,2k+2j-1)$ and from
$(1,2k+2j)$ to\break $(2k-1,2k+2j+1)$, for $1\le j<k$.

The final phase of the proof consists of establishing the relations
$$\openup1\jot
\eqalign{&\hbox{$j$ odd}\cr
\noalign{\smallskip}
(0,2j)&\sim(2k,2j-1)\cr
(1,2k+2j)&\sim\ol(2k,4k-2j-1)\cr
(2k,2j+1)&\sim(2k-1,2k+2j+1)\cr
(0,2j)&\sim(2k,2j+1)\cr}
\qquad\qquad
\eqalign{&\hbox{$j$ even}\cr
\noalign{\smallskip}
(0,2j)&\sim(1,2k+2j)\cr
(0,4k-2j)&\sim(2k,4k-2j-1)\cr
(0,4k-2j-2)&\sim\ol(2k-1,2k+2j+1)\cr
(0,2j)&\sim(2k,2j+1)\cr}$$
for $1\le j<k$. Suppose first that $j=1$ and $k>1$. We know already that
$(0,2)\sim(2k,1)$. Now $(2k,4k-3)$ cannot be joined to $(0,4k-4)$, because
that would make a short circuit; it cannot be joined to $(0,4k-2)$, because
the neighbors of $(0,4k-2)$ are known. So we have
$(2k,4k-3)\sim(2k-1,6k-3)=\ol(1,2k+2)$. This implies
$(1,2k+2)\not\sim(0,2)$, so $(0,2)\sim(2k,3)$. We also have a path
$$(2k-1,2k+3)\sim\cdots\sim(1,2k+2)\sim\ol(2k,4k-3)\sim\cdots\sim\ol(0,4k-4)\,;
$$
hence $(2k-1,2k+3)\not\sim(2k,4k+3)=\ol(0,4k-4)$. This forces $(2k-1,2k+3)
\sim(2k,3)$. The proof for $j=1$ is complete.

Suppose the relations have been proved for some $j<k-1$. If $j$ is odd, we
have $(0,4k-2j-2)\not\sim(2k,4k-2j-1)$ and $(0,4k-2j-2)\not\sim(1,6k-2j-2)=
\ol(2k-1,2k+2j+1)$; hence $(0,4k-2j-2)\sim(2k,4k-2j-3)$. Also
$(0,2j+2)\not\sim(2k,2j+1)$; we must have
$(2k,2j+3)\sim(0,2j+2)\sim(1,2k+2j+2)$. And there's a path
$$(2k-1,2k+2j+3)\sim\cdots\sim(1,2k+2j+2)\sim\cdots\sim(2k,2j+3)\,,$$
so $(2k-1,2k+2j+3)\not\sim(2k,2j+3)$. This forces
$(2k-1,2k+2j+3)\sim(2k,4k+2j+3)=\ol(0,4k-2j-4)$. The required relations for
$j+1$ have been established.

If $j$ is even, we have proved that $(2k,2j+1)\not\sim(2k-1,2k+2j+1)$, so
$(2k,2j+1)\sim(0,2j+2)$. Also $(2k,4k-2j-3)\not\sim(0,4k-2j-2)$, and
$(2k,4k-2j-3)\not\sim(0,4k-2j-4)$ because of a short circuit; so
$(2k,4k-2j-3)\sim(2k-1,6k-2j-3)=\ol(1,2k+2j+2)$. This implies
$(0,2j+2)\not\sim(1,2k+2j+2)$, so $(0,2j+2)\sim(2k,2j+3)$. Finally, the path
$$(2k-1,2k+2j+3)\sim\cdots\sim(1,2k+2j+2)\sim\ol(2k,4k-2j-3)\sim\cdots\sim
\ol(0,4k-2j-4)$$
shows that $(2k-1,2k+2j+3)\not\sim(2k,4k+2j+3)=\ol(0,4k-2j-4)$; we must have
$(2k-1,\allowbreak{2k+2j+3})\sim(2k,2j+3)$.

Now that the induction on $j$ is complete, we have deduced the entire
Hamiltonians circuit with the exception of one link from one vertex (and its
complement). More precisely,when $k>1$ is odd we have found the
Hamiltonian path
$$\openup1\jot
\eqalign{(2k,2k-1)\sim(0,2k-2)
 &\sim(1,4k-2)\sim\cdots\sim(2k-1,4k-1)\sim\ol(0,2k)\sim\cdots\sim
    \ol(2k,2k+1)\cr
 &\sim\ol(0,2k+2)\sim\cdots\sim\ol(2k,2k+3)\sim(1,4k-4)\sim\cdots\sim
    (2k-1,4k-3)\cr
 &\sim(2k,2k-3)\sim(0,2k-4)\sim(2k,2k-5)\cr
\sim\cdots&\sim(2k,1)\sim\cdots\sim\ol(2k,1)\sim\cdots\sim\ol(2k,2k-1)\,;\cr}$$
when $k$ is even we have found another,
$$\eqalign{(0,2k)\sim\cdots\sim(2k,2k+1)
 &\sim\ol(1,4k-2)\sim\cdots\sim\ol(2k-1,4k-1)\sim\ol(2k,2k-1)\sim\ol(0,2k-2)\cr
 &\sim\ol(2k,2k-3)\sim\ol(0,2k-4)\sim
   \ol(1,4k-4)\sim\cdots\sim\ol(2k-1,4k-3)\cr
 &\sim(0,2k+2)\sim\cdots\sim(2k,2k+3)\sim(0,2k+4)\cr
\sim\cdots&\sim(2k,4k-1)\sim\ol(2k,4k-1)\sim\cdots\sim\ol(0,2k)\,.\cr}$$
The endpoints are not adjacent, so it is impossible to complete a circuit.
\quad\pfbox

\medskip
Willcocks [\cite\jw] 
also conjectured that square boards with side $<2(r+s)$ do {\it not\/} yield
Hamiltonian graphs. 
Using a slight extension of the methods above we can in fact prove a bit more:

\proclaim
Theorem 8. An $\{r,s\}$-leaper has no Hamiltonian circuit on an $m\times n$
board when $2s\leq m\leq n<2(r+s)$.

\proof
We may assume that $2\leq r<s$. We show first that there is no Hamiltonian
circuit on an $m\times n$ board when $m=2s$ and $n$ is arbitrary. Say that
vertex $(x,y)$ is type~$A$ if $x<s$ and $x\equiv t\;(\bmod\;2r)$, or
if $x\geq s$ and $x\equiv r+s+t\;(\bmod \;2r)$, where
$$t=\cases{s\bmod 2r\,,&if $s\bmod 2r<r$;\cr
0\,&if $s\bmod 2r>r$.\cr}$$
Similarly, say that $(x,y)$ is type~$B$ if $x<s$ and $x\equiv r+t$,
 or $x\geq s$
and $x\equiv s+t\;(\bmod \;2r)$. Otherwise $(x,y)$ is type~$C$. Let $s=2kr+s'$,
where $0\leq s'<2r$. If $s'<r$, we have $t=s'$, so the vertices of type~$A$ are
those in ranks $t$, $t+2r$, $\ldots\,$, $t+(2k-2)r$, $2t+(2k+1)r$, $\ldots\,$,
$2t+(4k-1)r$, while those of type~$B$ are in ranks $t+r$, $t+3r$, $\ldots\,$,
$t+(2k-1)r$, $2t+2kr$, $\ldots\,$, $2t+(4k-2)r$. If $s'>r$, we have $t=0$, so
the vertices of type~$A$ have ranks~$0$, $2r$, $\ldots\,$, $2kr$, $s'+(2k+1)r$,
$\ldots\,$, $s'+(4k+1)r$ while those of type~$B$ have ranks $r$, $3r$,
$\ldots\,$, $(2k+1)r$, $s'+2kr$, $\ldots\,$, $s'+4kr$. In both cases there are
exactly as many vertices of type~$B$ as type~$A$, and every neighbor of a
type~$A$ vertex has type~$B$. This rules out a Hamiltonian circuit, as in
Theorem~7.

To complete the proof, we must show that no Hamiltonian circuit is possible on
an $m\times n$ board when $2s<m\leq n<2(r+s)$. Let $x=\min(m-2s,r)-1$,
$y=\min(n-2s,r)-1$. Then the short circuit
$$\displaylines{(x,y),\,(x+r,y+s),\,(x,y+2s),\,(x+s,y+2s-r),\,(x+2s,y+2s),\cr
(x+2s-r,y+s),\,(x+2s,y),\,(x+s,y+r),\,(x,y)\cr}$$
is forced. \ \pfbox

\medskip
It would be very interesting to see a proof of Willcocks's general Hamiltonian
circuit conjecture. 
 A~presumably simpler problem, but also of interest, is
to characterize the smallest boards on which leaper graphs are biconnected.
The next cases to consider are perhaps those in which $r$ and $s$ are
consecutive elements of the sequence 1, 2, 5, 12, 29, 70, 169, \dots\
defined by the recurrence $A_{k+1}=2A_k+A_{k-1}$. This choice makes the
angles between the eight leaper moves as nearly equal as possible;
we have proved the conjecture only in the cases where those angles are as
unequal as possible.

The results above suggest several additional open problems. For example,
what is the diameter of the $\{r,s\}$-leaper graph on an $(r+s)\times (2s)$
board, when the conditions of Theorem~1 are satisfied? What is the
3-dimensional analog of Theorem~1?

What is the smallest $n$ for which $\{1,2k\}$-leapers can make a Hamiltonian
circuit of a $(2k+1)\times n$ board? The proof of Theorem~7 shows that such
circuits have an intriguing structure. When $k=1$, the answer is 10 (see
Bergholt~[\cite\berg]), but for larger values of $k$ it appears likely that
the answer is~$12k$. This conjecture is true, at any rate, when $k=2$; also
$n\ge36$ is necessary when $k=3$.

Theorem 3 provides a lower bound for certain Hamiltonian graphs, but it is not
the best possible result of its kind. The lower bound on~$n$
can, for example, be
raised by~2 whenever $r$ is a multiple of~4, because we can extend the argument
in the proof as follows: Suppose $r=2k$ and $n=(s+4)s-1$. Then the files that
are not $R$ or~$L$ are $[a,a]$ for $1\leq a\leq k+2$; $[a+1,a]$ and $[a+2,a]$,
for $1\leq a\leq k+1$; $[a+3,a]$ for $1\leq a\leq k$; $[a,1]$ for $4<a\leq
k+2$; and the complements of those files. As before we call files $[a,b]$ and
$\overline{[a,b]}$ odd or even according as $a$ is odd or even. It turns out
that the odd files have $2r-4$ more unspecified neighbors than the even files
do. All links go between an odd file and an even file, except for a link from
$[k+2,1]$ to $\overline{[k+2,1]}$.
So the odd excess can be dissipated only if the special connection from
$[k+2,1]$ to $\overline{[k+2,1]}$ goes from odd to odd. A~Hamiltonian circuit
is therefore impossible if $k$ is even.

This result is best possible when $r=4$, because numerous $\{4,5\}$-leaper
tours exist on a $9\times 46$ board. Here, for example, is one that can be
found using the method 
Euler [\cite\euler]
proposed for ordinary knights:

\medskip

{
\newcount\n
\def\\{\afterassignment\olin\n=}
\def\olin{\expandafter\overline\expandafter{\number\n}}

\halign{&\hbox to1.9em{$\hfil#\hfil$}\cr
    1&\\133&\\ 60&\\103&\\ 94&\\117&  173&\\143&\\202&\\ 19&    9&\\ 41&\\ 76&\\ 85&  183&  156&   31&\\ 17&   11&\\ 27&\\190&\\ 49&   83&  162&  \ldots\cr
\\135&\\ 58&\\105&\\124&\\ 67&\\ 96&\\115&  175&\\141&\\206&\\ 43&\\ 78&\\ 87&\\126&\\ 69&  181&  158&\\204&  196&\\188&\\ 47&\\150&  164&  185&  \ldots\cr
\\ 56&\\107&\\122&  168&\\128&\\ 65&\\ 98&\\113&  177&\\137&\\ 80&\\ 89&\\148&  166&\\ 36&\\ 71&  179&\\139&\\186&\\ 45&  151&   46&  187&\\195&  \ldots\cr
\\109&\\120&  170&\\146&    4&\\130&\\ 63&\\100&\\111&\\ 54&\\ 91&  153&\\199&\\ 22&    6&\\ 38&\\ 73&\\ 52&\\163&\\ 82&   48&  189&\\197&\\ 24&  \ldots\cr
\\118&  172&\\144&\\201&\\ 20&    2&\\132&\\ 61&\\102&\\ 93&  155&   30&\\ 16&   12&\\ 26&    8&\\ 40&\\ 75&\\ 84&   50&  191&   28&\\ 14&   14&  \ldots\cr
\\ 95&\\116&  174&\\142&\\207&\\134&\\ 59&\\104&\\125&\\ 68&  182&  157&\\203&\\ 18&   10&\\ 42&\\ 77&\\ 86&  184&  161&   32&  193&   24&  197&  \ldots\cr
\\ 66&\\ 97&\\114&  176&\\136&\\ 57&\\106&\\123&  167&\\127&\\ 70&  180&\\140&\\205&\\ 44&\\ 79&\\ 88&\\149&  165&\\ 35&  159&   34&  195&\\187&  \ldots\cr
\\129&\\ 64&\\ 99&\\112&\\ 55&\\108&\\121&  169&\\147&    5&\\ 37&\\ 72&  178&\\138&\\ 81&\\ 90&  152&\\198&\\ 23&\\194&\\ 33&\\160&\\185&\\164&  \ldots\cr
    3&\\131&\\ 62&\\101&\\110&\\119&  171&\\145&\\200&\\ 21&    7&\\ 39&\\ 74&\\ 53&\\ 92&  154&   29&\\ 15&   13&\\ 25&\\192&\\ 51&\\162&\\ 83&  \ldots\cr}

}

\smallskip
\noindent
(Only the left portion of the board is shown; the right half is reversed and
complemented, so that the full tour has $180^{\circ}$ symmetry. A~bar over a
number means that 207 should be added.)

A similar argument shows that the lower bound for $r=5$ can be raised from 52
to~56, and that a symmetric $\{5,6\}$-leaper tour does exist on an $11\times
56$ board.

For $r\geq 6$, the lower bounds derived above are not optimum, but more
powerful methods will be needed to establish the best possible results.
Computer algorithms for the general symmetric
traveling salesrep problem show
that the $\{6,7\}$-leaper graph on a $13\times 76$
board is not Hamiltonian; in fact, at least
18~additional edges are
needed to make a Hamiltonian circuit possible. This result [\cite\jr]
was obtained and
verified by two independently developed computer codes, one by Giovanni Rinaldi
and Manfred Padberg, the other by Michael J\"unger, Gerhard Reinelt, and Stefan
Thienel. It is interesting to note that the lower bound of 18 can be obtained
simply from a linear programming model with the constraint that each vertex has
degree~2. If each edge of the graph has weight~0 and each edge not in the graph
has weight~1, a tour of total weight~18 turns out to be possible.

 Leaper graphs
should provide good challenges for all such computer codes. Michael
J\"unger [\cite\juengi] has recently used his program to verify Willcocks's
conjecture when $r+s=15$, so the smallest unsettled cases are now
$\{r,17-r\}$-leapers for $2\le r\le7$. J\"unger [\cite\juengii] has also found
the smallest Hamiltonian graphs of $\{r,r+1\}$-leapers on $(2r+1)\times n$
boards when $r=6$ ($n=92$) and $r=7$ ($n=106$).

\vfill\eject

\centerline{References}

\bib
[\cite\lucas]
\'Edouard Lucas, {\sl R\'ecr\'eations Math\'ematiques\/ \bf 4} (Paris, 1894).

\bib
[\cite\dickins]
Anthony Dickins, {\sl A Guide to Fairy Chess\/} (Richmond, Surrey: The Q Press,
1967).

\bib
[\cite\nashw]
C. St~J. A. Nash-Williams, ``Abelian groups, graphs and generalized
knights,'' {\sl Proceedings of the Cambridge Philosophical Society\/ \bf 55}
(1959), 232--238.

\bib
[\cite\jello]
George P. Jelliss, ``Theory of leapers,'' {\sl Chessics\/ \bf 2}, number~24
(Winter 1985), 86--98.

\bib
[\cite\jelli]
George P. Jelliss, ``Generalized knights and Hamiltonian tours,''
preprint dated 3~September 1993.

\bib
[\cite\dawson]
T. R. Dawson, problem in {\sl L'Echiquier\/} (Brussels), 1928, pages 985 and
1054.

\bib
[\cite\frolow]
M. Frolow, {\sl Les Carr\'es Magiques\/} (Paris, 1886).

\bib
[\cite\jw]
G. P. Jelliss and T. H. Willcocks, ``The five free leapers,'' {\sl Chessics\/
\bf 1}, number~2 (July 1976),~2; number~6 (August 1978), 4--5.

\bib
[\cite\euler]
L. Euler, ``Solution d'une question curieuse qui ne paroit soumise {\`a}
aucune analyse,'' {\sl M\'emoires de l'Academie Royale des Sciences et Belles
Letters\/} (Berlin, 1759), 310--337.

\bib
[\cite\hubers]
E. Huber-Stockar, problem 6304, {\sl Fairy Chess Review\/ \bf 5}, number~16
and~17 (1945), 124, 134.

\bib
[\cite\jaenisch]
C. F. de Jaenisch, {\sl Trait\'e des Applications de l'Analyse Math\'ematique
au Jeu des Echecs}, St.~Petersbourg, 1862--1863, volume~2.

\bib
[\cite\flye]
C. Flye Sainte-Marie, ``Note sur un probl\`eme relatif \`a la marche du
cavalier sur l'\'echiquier,'' {\sl Bulletin de la Soci\'et\'e Math\'ematique de
France\/ \bf 5} (1877), 144--150.

\bib
[\cite\berg]
E. Bergholt, Letter to the editor, {\sl The British Chess Magazine}, 1918,
page~74.

\bib
[\cite\jr]
Michael J\"unger and Giovanni Rinaldi, electronic communications from the
University of Co\-logne, 29~March 1994.

\bib
[\cite\juengi] Michael J\"unger, electronic communication from the University
of Cologne, 8~April 1994.

\bib
[\cite\juengii] Michael J\"unger, electronic communication from the University
of Cologne, 11~April 1994.

\bye